\documentclass{article}
\usepackage{amsmath}
\usepackage{amsthm}
\usepackage{mathrsfs}
\usepackage{amsfonts}
\begin {document}
\topmargin= -.2in \baselineskip=20pt

\newtheorem{theorem}{Theorem}[section]
\newtheorem{proposition}[theorem]{Proposition}
\newtheorem{lemma}[theorem]{Lemma}
\newtheorem{corollary}[theorem]{Corollary}

\theoremstyle{remark}
\newtheorem{remark}[theorem]{Remark}

\title {Functional Equations of
$L$-Functions for Symmetric Products of the Kloosterman Sheaf
\thanks{The research of Lei Fu is supported by the NSFC
(10525107).}}
\author {Lei Fu\\
{\small Institute of Mathematics, Nankai University, Tianjin, P.
R. China}\\
{\small leifu@nankai.edu.cn}\\{}\\
Daqing Wan\\
{\small Department of Mathematics, University of California, Irvine,
CA
92697}\\
{\small dwan@math.uci.edu}}

\date{}
\maketitle

\begin{abstract}
We determine the (arithmetic) local monodromy at $0$ and at $\infty$
of the Kloosterman sheaf using local Fourier transformations and
Laumon's stationary phase principle. We then calculate
$\epsilon$-factors for symmetric products of the Kloosterman sheaf.
Using Laumon's product formula, we get functional equations of
$L$-functions for these symmetric products, and prove a conjecture
of Evans on signs of constants of functional equations.



\end{abstract}

\section*{Introduction}
Let $p\not=2$ be a prime number and let $ {\mathbb F}_p$ be the
finite field with $p$ elements. Fix an algebraic closure ${\mathbb
F}$ of ${\mathbb F}_p$. Denote the projective line over ${\mathbb
F}_p$ by ${\mathbb P}^1$. For any power $q$ of $p$, let ${\mathbb
F}_q$ be the finite subfield of ${\mathbb F}$ with $q$ elements. Let
$\ell$ be a prime number different from $p$. Fix a nontrivial
additive character $\psi: {\mathbb F}_p\to \overline {\mathbb
Q}_\ell^\ast$. For any $x\in {\mathbb F}_q^\ast$, we define the one
variable Kloosterman sum by
$$\mathrm {Kl}_2({\mathbb F}_q, x)=\sum\limits_{x\in  {\mathbb F}_q^\ast}
\psi\left(\mathrm {Tr}_{{\mathbb F}_q/{\mathbb
F}_p}\left(\lambda+\frac{x}{\lambda}\right)\right).$$ In \cite{D1},
Deligne constructs a lisse $\overline{\mathbb Q}_l$-sheaf $\mathrm
{Kl}_2$ of rank $2$ on ${\mathbb G}_m={\mathbb P}^1-\{0, \infty\}$,
which we call the Kloosterman sheaf, such that for any $x\in
{\mathbb G}_m({\mathbb F}_q)={\mathbb F}_q^\ast$, we have
$$\mathrm {Tr}(F_x, \mathrm {Kl}_{2,\bar x})=-\mathrm {Kl}_2({\mathbb F}_q,
x),$$ where $F_x$ is the geometric Frobenius element at the point
$x$.  For a positive integer $k$, the $L$-function $L({\mathbb G}_m,
\mathrm {Sym}^k(\mathrm {Kl}_2),T)$ of the $k$-th symmetric product
of $\mathrm {Kl}_2$ was first studied by Robba \cite{R} via Dwork's
$p$-adic methods. Motivated by applications in coding theory, by
connections with modular forms, $p$-adic modular forms and Dwork's
unit root zeta functions, there has been a great deal of recent
interests to understand $L({\mathbb G}_m, \mathrm {Sym}^k(\mathrm
{Kl}_2),T)$ as much as possible for all $k$ and for all $p$. This
quickly raises a large number of interesting new problems.

Let $j:{\mathbb G}_m\to {\mathbb P}^1$ be the inclusion. We shall be
interested in the $L$-function
$$M_k(p,T):=L({\mathbb P}^1, j_\ast(\mathrm {Sym}^k(\mathrm {Kl}_2)),T).$$
This is the non-trivial factor of $L({\mathbb G}_m, \mathrm
{Sym}^k(\mathrm {Kl}_2),T)$. The trivial factor of $L({\mathbb G}_m,
\mathrm {Sym}^k(\mathrm {Kl}_2),T)$ was completely determined in
Fu-Wan \cite{FW1}. By general theory of Grothendieck-Deligne, the
non-trivial factor $M_k(p,T)$ is a polynomial in $T$ with integer
coefficients, pure of weight $k+1$. Its degree $\delta_k(p)$ can be
easily extracted from Fu-Wan \cite{FW} Proposition 2.3, Lemmas 4.1
and 4.2:
$$\delta_k(p) = \left\{\begin{array}{cl}
\frac{k-1}{2} -\left[\frac{k}{2p}+\frac{1}{2}\right] &\hbox{ if } k
\hbox { is odd},\\
2\left( [\frac{k-2}{4}] -[\frac{k}{2p}]\right)& \hbox { if } k
\hbox{ is even}.\end{array}\right.$$

For fixed $k$, the variation of $M_k(p,T)$ as $p$ varies should be
explained by an automorphic form, see Choi-Evans \cite{CE} and
Evans \cite{E} for the precise relations in the cases $k\leq 7$
and Fu-Wan \cite{FW2} for a motivic interpretation for all $k$.
For $k\leq 4$, the degree $\delta_k(p) \leq 1$ and $M_k(p,T)$ can
be determined easily. For $k=5$, the degree $\delta_5(p)=2$ for
$p>5$. The quadratic polynomial $M_5(p,T)$ is explained by an
explicit modular form \cite{PTV}. For $k=6$, the degree
$\delta_6(p)=2$ for $p>6$. The quadratic polynomial $M_6(p,T)$ is
again explained by an explicit modular form \cite{HS}. For $k=7$,
the degree $\delta_7(p)=3$ for $p>7$. The cubic polynomial
$M_7(p,T)$ is conjecturally explained in a more subtle way by an
explicit modular form in Evans \cite{E}. We will return to this
conjecture later in the introduction.

For fixed $p$, the variation of $M_k(p,T)$ as $k$ varies
$p$-adically should be related to $p$-adic automorphic forms and
$p$-adic $L$-functions. No progress has been made along this
direction. The $p$-adic limit of $M_k(p,T)$ as $k$ varies
$p$-adically links to an important example of Dwork's unit root zeta
function, see the introduction in Wan \cite{W}. The polynomial
$M_k(p,T)$ can be used to determine the weight distribution of
certain codes, see Moisio \cite{M1}\cite{M2}, and this has been
studied extensively for small $p$ and small $k$. The $p$-adic Newton
polygon (the $p$-adic slopes) of $M_k(p,T)$ remains largely
mysterious.

\bigskip
By Katz \cite{K} 4.1.11, we have $(\mathrm {Kl}_2)^\vee=\mathrm
{Kl}_2\otimes \overline {\mathbb Q}_\ell(1)$. So for any natural
number $k$, we have $$(\mathrm {Sym}^k(\mathrm {Kl}_2))^\vee=\mathrm
{Sym}^k(\mathrm {Kl}_2)\otimes \overline {\mathbb Q}_\ell(k).$$
General theory (confer \cite{L} 3.1.1) shows that $M_k(p,T)$
satisfies the functional equation
$$M_k(p,T)
=cT^\delta  M_k\left(p,\frac{1}{p^{k+1}T}\right),$$ where
\begin{eqnarray*}
c&=&\prod_{i=0}^2\mathrm {det}(-F,H^i({\mathbb P}^1_{\mathbb F},
j_\ast(\mathrm {Sym}^k(\mathrm {Kl}_2))
)^{(-1)^{i+1}},\\
\delta&=&-\chi({\mathbb P}^1_{\mathbb F}, j_\ast(\mathrm
{Sym}^k(\mathrm {Kl}_2)) =\delta_k(p),
\end{eqnarray*}
and $F$ denotes the Frobenius correspondence. Applying the
functional equation twice, we get
$$c^2=p^{(k+1)\delta}.$$ Based on numerical computation, Evans \cite{E} suggests
that the sign of $c$ should be $-\left(\frac{p}{105}\right)$ (the
Jacobi symbol) for $k=7$, and $-\left(\frac{p}{1155}\right)$ for
$k=11$. In this paper, we determine $c$ for all $k$ and all $p>2$.
The main result of this paper is the following theorem.

\begin{theorem} Let $p>2$ be an odd prime. If $k$ is even, we have
$$c = p^{(k+1)([\frac{k-2}{4}]-[\frac{k}{2p}])}.$$
If $k$ is odd, we have
$$c=(-1)^{\frac{k-1}{2}+[\frac{k}{2p}+\frac{1}{2}]}p^{\frac{k+1}{2}(\frac{k-1}{2}-[\frac{k}{2p}+\frac{1}{2}])}
\left( \frac{-2}{p}\right)^{[\frac{k}{2p}+\frac{1}{2}]}
\prod_{j\in\{0,1,\ldots,[\frac{k}{2}]\},\; p\not\; | 2j+1}\left(
\frac{(-1)^{j}(2j+1)}{p}\right).
$$
\end{theorem}

\begin{corollary} If $k$ is even and $p>2$, the sign of $c$ is always
$1$. If $k$ is odd and $p>k$, the sign of $c$ is
$$(-1)^{\frac{k-1}{2}}\prod_{j\in\{0,1,\ldots,[\frac{k}{2}]\},\; p\not\; | 2j+1}\left(
\frac{(-1)^{j}(2j+1)}{p}\right).$$
\end{corollary}

In the above corollary, if we take $k=7$, we see that the sign of
$c$ for $p>7$ is
$$-\left(\frac{1\cdot (-3)\cdot 5\cdot (-7)}{p}\right)=-\left(\frac{105}{p}\right)
=-\left(\frac{p}{105}\right);$$ if we take $k=11$, we see that the
sign of $c$ for $p>11$ is
$$-\left(\frac{1\cdot (-3)\cdot 5\cdot (-7)\cdot 9\cdot (-11)}{p}\right)
=-\left(\frac{-1155}{p}\right)=-\left(\frac{p}{1155}\right),$$
consistent with Evans' calculation.

In the case $k=7$, Evans proposed a precise description of
$M_7(p,T)$ in terms of modular forms. For $k=7$ and $p>7$, the
polynomial $M_7(p,T)$ has degree $3$. Write
$$M_7(p,T) = 1 + a_p T +d_pT^2 + e_pT^3.$$
 The
functional equation and our sign determination show that one of the
reciprocal roots for $M_7(p,T)$ is $\left(\frac{p}{105}\right)p^4$
and $e_p =-\left(\frac{p}{105}\right)p^{12}.$ Denote the other two
reciprocal roots by $\lambda_p$ and $\mu_p$ which are Weil numbers
of weight $8$. We deduce that
$$a_p =-\left(\left(\frac{p}{105}\right)p^4 +\lambda_p +\mu_p\right),
\ \lambda_p\mu_p =p^8, \ |\lambda_p|=|\mu_p|=p^4.$$ To explain the
numerical calculation of Evans, Katz suggests that there exists a
two dimensional representation  $$\rho:\mathrm {Gal}(\overline
{\mathbb Q}/{\mathbb Q})\to\mathrm {GL}(\overline{\mathbb
Q}_\ell^2)$$ unramified for $p>7$ and a Dirichlet character $\chi$
such that
\begin{eqnarray*}
\alpha_p^2&=&
\chi(p)\left(\frac{p}{105}\right)\frac{\lambda_p}{p^4},\\
\beta_p^2&=& \chi(p)\left(\frac{p}{105}\right)\frac{\mu_p}{p^4},\\
\alpha_p\beta_p&=&\chi(p),
\end{eqnarray*}
where $\alpha_p$ and $\beta_p$ are the eigenvalues of the geometric
Frobenius element at $p$ under $\rho$. We then have
\begin{eqnarray*}
1-\left(\frac{p}{105}\right)\frac{a_p}{p^4}&=&
2+\left(\frac{p}{105}\right)\frac{\lambda_p}{p^4}
+\left(\frac{p}{105}\right)\frac{\mu_p}{p^4}\\
&=&\overline \chi(p)(2\alpha_p\beta_p+\alpha_p^2+\beta_p^2)\\
&=&\overline\chi(p)(\alpha_p+\beta_p)^2.
\end{eqnarray*}
Set $b(p)=p(\alpha_p+\beta_p)$. Evans \cite{E} conjectured that
$b(p)$ is the $p$-th Hecke eigenvalue for a weight $3$ newform $f$
on $\Gamma_0(525)$. Our $a_p$ equals $-c_pp^2$ in \cite{E}.

\bigskip
Our proof of Theorem 0.1 naturally splits into two parts,
corresponding to the two ramification points at $0$ and $\infty$.
Let $t$ be the coordinate of ${\mathbb A}^1={\mathbb
P}^1-\{\infty\}$. For any closed point $x$ in ${\mathbb P}^1$, let
${\mathbb P}^1_{(x)}$ be the henselization of ${\mathbb P}^1$ at
$x$. By Laumon's product formula \cite{L} 3.2.1.1, we have
$$c=p^{k+1}\prod_{x\in |{\mathbb P}^1|} \epsilon(
{\mathbb P}^1_{(x)}, j_\ast(\mathrm {Sym}^k(\mathrm
{Kl}_2))|_{{\mathbb P}^1_{(x)}}, dt|_{{\mathbb P}^1_{(x)}}),$$
where $|{\mathbb P}^1|$ is the set of all closed points of
${\mathbb P}^1$. When $x\not=0,\infty$, the sheaf $\mathrm
{Sym}^k(\mathrm {Kl}_2)|_{{\mathbb P}^1_{(x)}}$ is lisse and the
order of $dt$ at $x$ is $0$. So by \cite{L} 3.1.5.4 (ii) and (v),
we have $$\epsilon( {\mathbb P}^1_{(x)}, j_\ast(\mathrm
{Sym}^k(\mathrm {Kl}_2))|_{{\mathbb P}^1_{(x)}}, dt|_{{\mathbb
P}^1_{(x)}})=1$$ for $x\not=0,\infty$. Therefore
$$c=p^{k+1} \epsilon(
{\mathbb P}^1_{(0)}, j_\ast(\mathrm {Sym}^k(\mathrm
{Kl}_2))|_{{\mathbb P}^1_{(0)}}, dt|_{{\mathbb P}^1_{(0)}})\epsilon(
{\mathbb P}^1_{(\infty)}, j_\ast(\mathrm {Sym}^k(\mathrm
{Kl}_2))|_{{\mathbb P}^1_{(\infty)}}, dt|_{{\mathbb
P}^1_{(\infty)}}).$$

In \S 1, we prove the following.

\begin{proposition} We have $$\epsilon(
{\mathbb P}^1_{(0)}, j_\ast(\mathrm {Sym}^k(\mathrm
{Kl}_2))|_{{\mathbb P}^1_{(0)}}, dt|_{{\mathbb
P}^1_{(0)}})=(-1)^kp^{\frac{k(k+1)}{2}}.$$
\end{proposition}

In \S 2, we prove the following.

\begin{proposition}
$\epsilon( {\mathbb P}^1_{(\infty)}, j_\ast(\mathrm {Sym}^k(\mathrm
{Kl}_2))|_{{\mathbb P}^1_{(\infty)}}, dt|_{{\mathbb
P}^1_{(\infty)}})$ equals
$$p^{-(k+1)(\frac{k+8}{4}+[\frac{k}{2p}])}$$
if $k=2r$ for an even $r$,
$$p^{-(k+1)(\frac{k+6}{4}+[\frac{k}{2p}])}$$
if $k=2r$ for an odd $r$, and
$$(-1)^{\frac{k+1}{2}+[\frac{k}{p}]-[\frac{k}{2p}]}p^{-\frac{k+1}{2}(\frac{k+5}{2}+[\frac{k}{p}]-[\frac{k}{2p}])}
\left( \frac{-2}{p}\right)^{[\frac{k}{p}]-[\frac{k}{2p}]}
\prod_{j\in\{0,1,\ldots,[\frac{k}{2}]\},\; p\not\; | 2j+1}\left(
\frac{(-1)^{j}(2j+1)}{p}\right)
$$
if $k=2r+1$.
\end{proposition}

We deduce from the above two propositions the constant $c$ as stated
in Theorem 0.1 using the following facts:
\begin{eqnarray*}
&&\left[\frac{k-2}{4}\right]= \left\{
\begin{array}{cl}
\frac{k-4}{4} &\hbox { if } k=2r \hbox { for an even }r,\\
\frac{k-2}{4} &\hbox { if } k=2r \hbox { for an odd }r,
\end{array}\right.\\
&&\left[\frac{k}{p}\right]-\left[\frac{k}{2p}\right]
=\left[\frac{k}{2p}+\frac{1}{2}\right] \hbox { if } k \hbox { is
odd}.
\end{eqnarray*}

To get Proposition 0.4, we first have to determine the local
(arithmetic) monodromy of $\mathrm {Kl}_2$ at $\infty$. This is
Theorem 2.1 in \S 2, which is of interest itself, and is proved by
using local Fourier transformations and Laumon's stationary phase
principle.

\section{Calculation of $\epsilon( {\mathbb
P}^1_{(0)}, j_\ast(\mathrm {Sym}^k(\mathrm {Kl}_2))|_{{\mathbb
P}^1_{(0)}}, dt|_{{\mathbb P}^1_{(0)}})$}

Let $\eta_0$ be the generic point of ${\mathbb P}^1_{(0)}$, let
$\bar\eta_0$ be a geometric point located at $\eta_0$, and let $V$
be an $\overline {\mathbb Q}_\ell$-representation of $\mathrm
{Gal}(\bar\eta_0/\eta_0)$. Suppose the inertia subgroup $I_0$ of
$\mathrm {Gal}(\bar\eta_0/\eta_0)$ acts unipotently on $V$. Consider
the $\ell$-adic part of the cyclotomic character
$$t_\ell: I_0\to {\mathbb Z}_\ell(1), \; \sigma\mapsto
\left(\frac{\sigma(\sqrt[\ell^n]{t})}{\sqrt[\ell^n]{t}}\right).$$
Note that for any $\sigma$ in the inertia subgroup, the $\ell^n$-th
root of unity $\frac{\sigma(\sqrt[\ell^n]{t})}{\sqrt[\ell^n]{t}}$
does not depends on the choice of the $\ell^n$-th root
$\sqrt[\ell^n]{t}$ of $t$. Since $I_0$ acts on $V$ unipotently,
there exists a nilpotent homomorphism
$$N:V(1)\to V$$ such that the action of $\sigma\in I_0$ on $V$
is given by $\exp (t_\ell(\sigma).N)$. Fix a lifting $F\in \mathrm
{Gal}(\bar\eta_0/\eta_0)$ of the geometric Frobenius element in
$\mathrm {Gal}({\mathbb F}/{\mathbb F}_p)$.

\begin{lemma}
Notation as above. Let $V=\mathrm {Kl}_{2,\bar\eta_0}$. There exists
a basis $\{e_0,e_1\}$ of $V$ such that
$$F(e_0)=e_0,\; F(e_1)=pe_1$$
$$N(e_0)=0,\; N(e_1)=e_0.$$
\end{lemma}

\begin{proof} This is the $n=2$ case of Proposition 1.1 in
\cite{FW1}.
\end{proof}

\begin{lemma} Keep the notation in Lemma 1.1. Let $\{f_0,\ldots,
f_k\}$ be the basis of $\mathrm {Sym}^k(V)=\mathrm {Sym}^k(\mathrm
{Kl}_{2,\bar\eta_0})$ defined by $f_i=\frac{1}{i!}e_0^{k-i}e_1^i.$
We have
$$F(f_i)=p^i f_i, N(f_i)=f_{i-1},$$
where we regard $f_{i-1}$ as $0$ if $i=0$.

\end{lemma}
\begin{proof} Use the fact that for any $v_1,\ldots, v_k\in V$,
we have the following identities in $\mathrm {Sym}^k(V)$:
\begin{eqnarray*}
F(v_1\cdots v_k)&=&F(v_1)\cdots F(v_k), \\
N(v_1\cdots v_k)&=&\sum_{i=1}^k v_1\cdots v_{i-1} N(v_i)
v_{i+1}\cdots v_k.
\end{eqnarray*}
\end{proof}

\begin{corollary} The sheaf $\mathrm {Sym}^k(\mathrm
{Kl}_2)|_{\eta_0}$ has a filtration $$0={\mathcal
F}_{-1}\subset{\mathcal F}_0\subset\cdots\subset {\mathcal
F}_k=\mathrm {Sym}^k(\mathrm {Kl}_2)|_{\eta_0}$$ such that
$${\mathcal F}_i/{\mathcal F}_{i-1}\cong \overline {\mathbb
Q}_l(-i)$$ for any $i=0,\ldots,k$.
\end{corollary}

\begin{proof} This follows from Lemma 1.2 by taking ${\mathcal F}_i$
to be the sheaf on $\eta_0$ corresponding to the galois
representation $\mathrm {Span}(f_0,\ldots, f_i)$ of $\mathrm
{Gal}(\bar\eta_0/\eta_0)$.
\end{proof}

The following is Proposition 0.3 in the introduction.

\begin{proposition} We have $$\epsilon(
{\mathbb P}^1_{(0)}, j_\ast(\mathrm {Sym}^k(\mathrm
{Kl}_2))|_{{\mathbb P}^1_{(0)}}, dt|_{{\mathbb
P}^1_{(0)}})=(-1)^kp^{\frac{k(k+1)}{2}}.$$
\end{proposition}

\begin{proof} Let $u:\eta_0\to {\mathbb P}^1_{(0)}$ and $v:
\{0\}\to {\mathbb P}^1_{(0)}$ be the immersions. By \cite{L} 3.1.5.4
(iii) and (v), we have
\begin{eqnarray*}
\epsilon( {\mathbb P}^1_{(0)}, u_\ast\overline {\mathbb Q}_\ell(-i),
dt|_{{\mathbb P}^1_{(0)}})&=&1,\\
\epsilon( {\mathbb P}^1_{(0)}, v_\ast\overline {\mathbb Q}_\ell(-i),
dt|_{{\mathbb P}^1_{(0)}})&=&\mathrm {det}(-F_0,\overline {\mathbb
Q}_\ell(-i))^{-1}=-\frac{1}{p^i}.
\end{eqnarray*}
We have an exact sequence
$$0\to u_! \overline
{\mathbb Q}_\ell(-i) \to u_\ast\overline {\mathbb Q}_\ell(-i)\to
v_\ast\overline {\mathbb Q}_\ell(-i)\to 0.$$ It follows from
\cite{L} 3.1.5.4 (ii) that we have
\begin{eqnarray*}
\epsilon( {\mathbb P}^1_{(0)}, u_!\overline {\mathbb Q}_\ell(-i),
dt|_{{\mathbb P}^1_{(0)}})&=& \frac{\epsilon( {\mathbb P}^1_{(0)},
u_\ast\overline {\mathbb Q}_\ell(-i), dt|_{{\mathbb P}^1_{(0)}})}{
\epsilon( {\mathbb P}^1_{(0)}, v_\ast\overline {\mathbb Q}_\ell(-i),
dt|_{{\mathbb P}^1_{(0)}})}\\
&=& -p^i.
\end{eqnarray*}
By Corollary 1.3 and \cite {L} 3.1.5.4 (ii), we have
\begin{eqnarray*}
\epsilon({\mathbb P}^1_{(0)}, j_!(\mathrm {Sym}^k(\mathrm
{Kl}_2))|_{{\mathbb P}^1_{(0)}}, dt|_{{\mathbb
P}^1_{(0)}})&=&\prod_{i=0}^k \epsilon( {\mathbb P}^1_{(0)},
u_!({\mathcal F}_i/{\mathcal F}_{i+1}), dt|_{{\mathbb P}^1_{(0)}})\\
&=& \prod_{i=0}^k \epsilon( {\mathbb P}^1_{(0)}, u_!{\mathbb
Q}_\ell(-i), dt|_{{\mathbb P}^1_{(0)}})\\
&=&\prod_{i=0}^k (-p^i).
\end{eqnarray*}
Moreover, by Lemma 1.2, we have $$v^\ast(j_\ast(\mathrm
{Sym}^k(\mathrm {Kl}_2))|_{{\mathbb P}^1_{(0)}})\cong \overline
{\mathbb Q}_\ell,$$ and hence
$$\epsilon( {\mathbb P}^1_{(0)}, v_\ast
v^\ast(j_\ast(\mathrm {Sym}^k(\mathrm {Kl}_2))|_{{\mathbb
P}^1_{(0)}}), dt|_{{\mathbb P}^1_{(0)}})=-1.$$ So we have
\begin{eqnarray*}
&&\epsilon( {\mathbb P}^1_{(0)}, j_\ast(\mathrm {Sym}^k(\mathrm
{Kl}_2))|_{{\mathbb P}^1_{(0)}}, dt|_{{\mathbb
P}^1_{(0)}})\\&=&\epsilon( {\mathbb P}^1_{(0)}, j_!(\mathrm
{Sym}^k(\mathrm {Kl}_2))|_{{\mathbb P}^1_{(0)}}, dt|_{{\mathbb
P}^1_{(0)}})\epsilon( {\mathbb P}^1_{(0)}, v_\ast
v^\ast(j_\ast(\mathrm {Sym}^k(\mathrm {Kl}_2))|_{{\mathbb
P}^1_{(0)}}), dt|_{{\mathbb P}^1_{(0)}})\\
&=&\prod_{i=1}^k (-p^i)\\
&=&(-1)^kp^{\frac{k(k+1)}{2}}.
\end{eqnarray*}
\end{proof}

\section{Calculation of $\epsilon( {\mathbb
P}^1_{(\infty)}, j_\ast(\mathrm {Sym}^k(\mathrm {Kl}_2))|_{{\mathbb
P}^1_{(\infty)}}, dt|_{{\mathbb P}^1_{(\infty)}})$}

We first introduce some notations. Fix a nontrivial additive
character $\psi:{\mathbb F}_p\to \overline{\mathbb Q}_\ell^\ast$ and
define $\mathrm {Kl}_2$ as in the introduction. Fix a separable
closure $\overline{{\mathbb F}_p(t)}$ of ${\mathbb F}_p(t)$. Let $x$
be an element in $\overline {{\mathbb F}_p(t)}$ satisfying
$x^p-x=t$. Then ${\mathbb F}_p(t,x)$ is galois over ${\mathbb
F}_p(t)$. We have a canonical isomorphism
$${\mathbb F}_p\stackrel {\cong}\to \mathrm {Gal}({\mathbb
F}_p(t,x)/{\mathbb F}_p(t))$$ which sends each $a\in {\mathbb F}_p$
to the element in $\mathrm {Gal}({\mathbb F}_p(t,x)/{\mathbb
F}_p(t))$ defined by $x\mapsto x+a$. Let ${\mathcal L}_\psi$ be the
galois representation defined by
$$\mathrm {Gal}(\overline{{\mathbb F}_p(t)}/{\mathbb F}_p(t))\to
\mathrm {Gal}({\mathbb F}_p(t,x)/{\mathbb F}_p(t))\stackrel
{\cong}\to {\mathbb F}_p \stackrel {\psi^{-1}}\to \overline {\mathbb
Q}_\ell^\ast.$$ It is unramfied outside $\infty$ and totally wild at
$\infty$ with Swan conductor $1$. This galois representation defines
a lisse $\overline {\mathbb Q}_\ell$-sheaf on ${\mathbb A}^1$ which
we still denote by ${\cal L}_{\psi}$. Let $X$ be an ${\mathbb
F}_p$-scheme. Any section $f$ in ${\mathcal O}_X(X)$ defines an
${\mathbb F}_p$-algebra homomorphism $${\mathbb F}_p[t]\to {\mathcal
O}_X(X),\; t\mapsto f,$$ and hence an ${\mathbb F}_p$-morphism of
schemes $$f:X\to {\mathbb A}^1.$$ We denote the lisse
$\overline{\mathbb Q}_\ell$-sheaf $f^\ast{\mathcal L}_\psi$ on $X$
by ${\mathcal L}_\psi(f)$. For any $f_1,f_2\in {\mathcal O}_X(X)$,
we have $${\mathcal L}_\psi(f_1)\otimes {\mathcal L}_\psi(f_2)\cong
{\mathcal L}_\psi(f_1+f_2).$$

Recall that $p\not =2$. Let $y$ be an element in $\overline
{{\mathbb F}_p(t)}$ satisfying $y^2=t$. Then ${\mathbb F}_p(t,y)$ is
galois over ${\mathbb F}_p(t)$. We have a canonical isomorphism
$$\{\pm 1\} \stackrel {\cong}\to \mathrm {Gal}({\mathbb
F}_p(t,y)/{\mathbb F}_p(t))$$ which sends $-1$ to the element in
$\mathrm {Gal}({\mathbb F}_p(t,y)/{\mathbb F}_p(t))$ defined by
$y\mapsto -y$. Let $$\chi:\{\pm 1\}\to \overline {\mathbb
Q}_\ell^\ast$$ be the (unique) nontrivial character. Define
${\mathcal L}_\chi$ to be the galois representation defined by
$$\mathrm {Gal}(\overline{{\mathbb F}_p(t)}/{\mathbb F}_p(t))
\to\mathrm {Gal}({\mathbb F}_p(t,y)/{\mathbb F}_p(t))\stackrel
{\cong}\to \{\pm 1\}\stackrel {\chi^{-1}}\to \overline {\mathbb
Q}_\ell^\ast.$$ It is unramified outside $0$ and $\infty$, and
tamely ramified at $0$ and $\infty$. This galois representation
defines a lisse $\overline {\mathbb Q}_\ell$-sheaf on ${\mathbb
G}_m$ which we still denote by ${\mathcal L}_\chi$.

Let $\theta:\mathrm {Gal}({\mathbb F}/{\mathbb F}_p)\to \overline
{\mathbb Q}_\ell^\ast$ be a character of the galois group of the
finite field. Denote by ${\mathcal L}_\theta$ the galois
representation
$$\mathrm {Gal}(\overline {{\mathbb F}_p(t)}/{\mathbb F}_p(t))
\to \mathrm {Gal}({\mathbb F}/{\mathbb F}_p) \stackrel {\theta} \to
\overline {\mathbb Q}_\ell^\ast.$$  It is unramified everywhere, and
hence defines a lisse $\overline {\mathbb Q}_l$-sheaf on ${\mathbb
P}^1$ which we still denote by ${\mathcal L}_\theta$.

\begin{theorem} Notation as above. Let $\eta_\infty$ be the generic
point of ${\mathbb P}^1_{(\infty)}$. Then $\mathrm
{Kl}_2|_{\eta_\infty}$ is isomorphic to the restriction to
$\eta_\infty$ of the sheaf
$$[2]_\ast({\cal L}_\psi(2t)\otimes {\cal L}_\chi) \otimes
{\cal L}_{\theta_0},$$ where $[2]:{\mathbb G}_m\to {\mathbb G}_m$ is
the morphism defined by $x\mapsto x^2$, and
$$\theta_0:\mathrm {Gal}({\mathbb F}/{\mathbb F}_p)\to
\overline{\mathbb Q}_\ell^\ast$$ is the character sending the
geometric Frobenius element $F$ in $\mathrm {Gal}({\mathbb
F}/{\mathbb F}_p)$ to the Gauss sum
$$\theta_0(F)=g(\chi,\psi)=-\sum_{x\in {\mathbb
F}_p^\ast}\left(\frac{x}{p}\right)\psi(x).$$
\end{theorem}

\begin{proof} By \cite{FW2} Proposition 1.1, we have
\begin{eqnarray}
\mathrm {Kl}_2={\mathcal F}\left(j_!{\mathcal L}_\psi
\left(\frac{1}{t}\right)\right)|_{{\mathbb G}_m},
\end{eqnarray}
where ${\mathcal F}$ is the $\ell$-adic Fourier transformation and
$j:{\mathbb G}_m\to {\mathbb A}^1$ is the inclusion. Let $$\pi_1,
\pi_2:{\mathbb G}_{m}\times_{{\mathbb F}_p} {\mathbb G}_{m}\to
{\mathbb G}_{m}$$ be the projections. Using the proper base change
theorem and the projection formula (\cite{SGA 4} XVII 5.2.6 and
5.2.9), one can verify
\begin{eqnarray}
[2]^\ast \left({\cal F}\left(j_!{\cal
L}_\psi\left(\frac{1}{t}\right)\right)|_{{\mathbb
G}_{m}}\right)\cong R\pi_{2!} \left({\cal
L}_\psi\left(\frac{1}{t}+tt'^{2}\right)\right)[1],
\end{eqnarray}
where $$\frac{1}{t}+tt'^{2}:{\mathbb G}_{m}\times_{{\mathbb
F}_p}{\mathbb G}_{m}\to {\mathbb A}^1$$ is the morphism
corresponding to the ${\mathbb F}_p$-algebra homomorphism
$${\mathbb F}_p[t]\to{\mathbb F}_p[t,1/t, t', 1/t'],\; t\mapsto \frac{1}{t}
+tt'^{2}.$$ Consider the isomorphism
$$\tau:{\mathbb G}_{m}\times_{{\mathbb F}_p}
{\mathbb G}_{m}\to {\mathbb G}_{m}\times_{{\mathbb F}_p} {\mathbb
G}_{m},\; (t,t')\mapsto \left(\frac{t}{t'},t'\right).$$ We have
$\pi_2\tau=\pi_2$. So
\begin{eqnarray}
R\pi_{2!} \left({\mathcal
L}_\psi\left(\frac{1}{t}+tt'^{2}\right)\right)\cong
R(\pi_2\tau)_!\tau^\ast \left({\mathcal
L}_\psi\left(\frac{1}{t}+tt'^{2}\right)\right) \cong
R\pi_{2!}{\mathcal L}_\psi\left(\left(\frac{1}{t}+t\right)t'\right).
\end{eqnarray}
Consider the morphism
$$g:{\mathbb G}_m\to {\mathbb A}^1, \; t\mapsto \frac{1}{t}+t.$$
Again using the proper base change theorem and the projection
formula, one can verify
\begin{eqnarray}
{\mathcal F}(Rg_!\overline {\mathbb Q}_\ell)\cong R\pi_{2!}{\cal
L}_\psi\left(\left(\frac{1}{t}+t\right)t'\right)[1].
\end{eqnarray}
From the isomorphisms (1)-(4), we get
$$
[2]^\ast \mathrm {Kl}_2 \cong {\mathcal F}(Rg_!\overline {\mathbb
Q}_\ell)|_{{\mathbb G}_m}.$$ By Lemma 2.2 below, the stationary
phase principle of Laumon \cite{L} 2.3.3.1 (iii), and \cite{L}
2.5.3.1, as representations of ${\rm
Gal}(\bar\eta_{\infty'}/\eta_{\infty'})$, we have
\begin{eqnarray*}
{\mathcal H}^0({\mathcal F}(Rg_!\overline{\mathbb Q}_\ell))_{\bar
\eta_{\infty'}} &\cong& {\mathcal F}^{(2,\infty')}({\mathcal
L}_{\chi})
\bigoplus {\mathcal F}^{(-2,\infty')}({\mathcal L}_{\chi})\\
&\cong& ({\mathcal L}_\psi(2t')\otimes  {\cal F}^{(0,\infty')}({\cal
L}_{\chi})) \bigoplus ({\mathcal L}_\psi(-2t')\otimes  {\cal
F}^{(0,\infty')}({\mathcal L}_{\chi}))
\\
&\cong& ({\mathcal L}_\psi(2t')\otimes {\mathcal L}_\chi\otimes
{\mathcal L}_{\theta_0})\bigoplus ({\mathcal L}_\psi(-2t')\otimes
{\mathcal L}_\chi\otimes {\mathcal L}_{\theta_0}).
\end{eqnarray*}
Hence
$$([2]^\ast \mathrm {Kl}_2)|_{\eta_\infty}\cong ({\mathcal L}_\psi(2t)\otimes {\mathcal L}_\chi\otimes
{\mathcal L}_{\theta_0})|_{\eta_\infty}\bigoplus ({\mathcal
L}_\psi(-2t)\otimes {\mathcal L}_\chi\otimes {\mathcal
L}_{\theta_0})|_{\eta_\infty}.$$  Note that this decomposition of
$([2]^\ast \mathrm {Kl}_2)|_{\eta_\infty}$ is non-isotypical. By
\cite{Se} Proposition 24 on p. 61, and the fact that $\mathrm
{Kl}_2|_{\eta_\infty}$ is irreducible (since its Swan conductor is
$1$), we have $$\mathrm {Kl}_2|_{\eta_\infty}\cong [2]_\ast
({\mathcal L}_\psi(2t)\otimes {\mathcal L}_\chi\otimes {\mathcal
L}_{\theta_0})|_{\eta_\infty}.$$ We have
\begin{eqnarray*}
[2]_\ast ({\mathcal L}_\psi(2t)\otimes {\mathcal L}_\chi\otimes
{\mathcal L}_{\theta_0}) &\cong&[2]_\ast\left({\mathcal
L}_\psi(2t)\otimes {\mathcal L}_\chi
\otimes [2]^\ast {\mathcal L}_{\theta_0}\right)\\
&\cong&  [2]_\ast\left({\mathcal L}_\psi(2t)\otimes {\mathcal
L}_\chi\right) \otimes {\mathcal L}_{\theta_0}.
\end{eqnarray*} Here we
use the fact that $[2]^\ast {\mathcal L}_{\theta_0}\cong {\mathcal
L}_{\theta_0}.$ Hence $$\mathrm {Kl}_2|_{\eta_\infty}\cong
\biggl([2]_\ast({\cal L}_\psi(2t)\otimes {\cal L}_\chi) \otimes
{\cal L}_{\theta_0}\biggr)|_{\eta_\infty}.$$
\end{proof}

\begin{lemma} For the morphism $$g:{\mathbb G}_m\to{\mathbb A}^1,\;
t\mapsto \frac{1}{t}+t,$$  the following holds:

(i) $Rg_!\overline{\mathbb Q}_\ell$ is a $\overline{\mathbb
Q}_\ell$-sheaf on ${\mathbb A}^1$ which is lisse outside the
rational points $2$ and $-2$.

(ii)  $Rg_!\overline{\mathbb Q}_\ell$ is unramified at $\infty$.

(iii) Let $P$ be one of the rational points $2$ or $-2$, and let
$\tilde {\mathbb A}^1_{(P)}$ be the henselization of ${\mathbb A}^1$
at $P$.  We have
$$(Rg_!\overline{\mathbb Q}_\ell)|_{\tilde {\mathbb A}^1_{(P)}}\cong
\overline{\mathbb Q}_\ell\oplus {\mathcal L}_{\chi, !},$$ where
${\mathcal L}_{\chi, !}$ denotes the extension by $0$ of the Kummer
sheaf ${\mathcal L}_{\chi}$ on the generic point of  $\tilde
{\mathbb A}^1_{(P)}$ to $\tilde {\mathbb A}^1_{(P)}$ .
\end{lemma}

\begin{proof} We have
$$\frac{\partial g}{\partial t}=-\frac{1}{t^2}+1.$$
So $\frac{\partial g}{\partial t}$ vanishes at the points $t=\pm 1$.
We have
\begin{eqnarray*}
&&g(\pm 1)=\pm 2, \\
&&\frac {\partial^2 g}{\partial t^2}(\pm 1)=\pm 2\not=0.
\end{eqnarray*}
It follows that $g$ is tamely ramified above $\pm 2$ with
ramification index $2$, and $g$ is \'etale elsewhere. Consider the
morphism
$$\bar g:{\mathbb P}^1\to{\mathbb P}^1,\; [t_0:t_1]\mapsto [t_0t_1:
t_0^2+t_1^2].$$ We have $\bar g^{-1}(\infty)=\{0,\infty\}.$ Hence
$$\bar g^{-1}({\mathbb A}^1)={\mathbb G}_{m}.$$
It is clear that $$\bar g|_{{\mathbb G}_{m}}=g.$$ So $g:{\mathbb
G}_{m}\to{\mathbb A}^1$ is a finite morphism of degree $2$. Near
$0$, the morphism $\bar g$ can be expressed as
$$t\mapsto \frac{t}{1+t^{2}}.$$ Hence $\bar g$ is unramified at $0$.
Similarly $\bar g$ is also unramified at $\infty$. Our lemma follows
from these facts.
\end{proof}

\begin{remark} The first attempt to determine the monodromy at
$\infty$ of the $(n-1)$-variable Kloosterman sheaf $\mathrm
{Kl}_n|_{\eta_\infty}$ is done in Fu-Wan \cite{FW} Theorem 1.1,
where we deduce from Katz \cite{K} that
$$\mathrm {Kl}_n|_{\eta_\infty}\cong \left([n]_\ast({\mathcal L}_\psi(nt)\otimes
{\mathcal L}_{\chi^{n-1}})\otimes {\mathcal
L}_{\theta}\otimes\overline{\mathbb
Q}_\ell\left(\frac{1-n}{2}\right)\right)|_{\eta_\infty}$$ for some
character $\theta: \mathrm {Gal}({\mathbb F}/{\mathbb F}_p)\to
\overline {\mathbb Q}_\ell^\ast$, and an explicit description of
$\theta^2$ is given. Using induction on $n$, \cite{FW2} Proposition
1.1, and adapting the argument in \cite{F} to non-algebraically
closed ground field, we can get an explicit description of $\theta$.
See \cite{F} where the monodromy of the more general hypergeometric
sheaf is treated (over algebraically closed field).
\end{remark}

\begin{lemma} Keep the notation in Theorem 2.1. Let
$$\theta_1: \mathrm {Gal}({\mathbb F}/{\mathbb F}_p)\to
\overline{\mathbb Q}_\ell^\ast$$ be the character defined by
$$\theta_1(\sigma)=\chi\left(\frac{\sigma(\sqrt{-1})}{\sqrt{-1}}\right)$$
for any $\sigma \in \mathrm {Gal}({\mathbb F}/{\mathbb F}_p)$. Note
that the above expression is independent of the choice of the square
root $\sqrt{-1}\in {\mathbb F}$ of ${-1}$.

(i) If $k=2r$ is even, $\mathrm {Sym}^k(\mathrm
{Kl}_2)|_{\eta_\infty}$ is isomorphic to the restriction to
$\eta_\infty$ of the sheaf
$$\left ({\mathcal L}_{{\chi}^r} \otimes
{\mathcal L}_{\theta_0^{2r}\theta_1^r} \right)\oplus
\left(\bigoplus_{i=0}^{r-1} [2]_\ast{\mathcal
L}_\psi((4i-4r)t)\otimes {\mathcal L}_{\theta_0^{2r}
\theta_1^{i}}\right).
$$

(ii) If $k=2r+1$ is odd, $\mathrm {Sym}^k(\mathrm
{Kl}_2)|_{\eta_\infty}$ is isomorphic to the restriction to
$\eta_\infty$ of the sheaf
$$\bigoplus_{i=0}^{r} [2]_\ast\left({\mathcal L}_\psi((4i-4r-2)t)\otimes {\mathcal
L}_{\chi}\right) \otimes {\mathcal
L}_{\theta_0^{2r+1}\theta_1^{i+1}}.$$
\end{lemma}

\begin{proof} By Theorem 2.1, it suffices to calculate the restriction to
$\eta_\infty$ of $\mathrm {Sym}^k([2]_\ast ({\mathcal
L}_{\psi}(2t)\otimes {\mathcal L}_{\chi}))$. Let $y$, $z$, $w$ be
elements in $\overline {{\mathbb F}_p(t)}$ satisfying $$y^2=t,\;
z^p-z=y,\; w^2=y.$$ Fix a square root $\sqrt{-1}$ of $-1$ in
${\mathbb F}$. Then ${\mathbb F}_p(z, w,\sqrt{-1})$ and ${\mathbb
F}_p(y)$ are galois extensions of ${\mathbb F}_p(t)$. Let $G=\mathrm
{Gal}({\mathbb F}_p(z,w,\sqrt{-1})/{\mathbb F}_p(t))$ and $H=\mathrm
{Gal}({\mathbb F}_p(z,w,\sqrt{-1})/{\mathbb F}_p(y))$. Then $H$ is
normal in $G$, and we have canonical isomorphisms
$$G/H\stackrel\cong\to \mathrm {Gal}({\mathbb F}_p(y)/{\mathbb F}_p(t))
\stackrel\cong \to \{\pm 1\}.$$ Consider the case where $\sqrt{-1}$
does not lie in ${\mathbb F}_p$. We have an isomorphism
$${\mathbb F}_p \times \{\pm 1\}\times \{\pm 1\} \stackrel \cong\to
H=\mathrm {Gal}({\mathbb F}_p(z,w,\sqrt {-1})/{\mathbb F}_p(y))$$
which maps $(a,\mu',\mu'')\in {\mathbb F}_p\times \{\pm
1\}\times\{\pm 1\}$ to the element $g_{(a,\mu',\mu'')}\in \mathrm
{Gal}({\mathbb F}_p(z,w,\sqrt{-1})/{\mathbb F}_p(y))$ defined by
$$g_{(a,\mu',\mu'')}(z)=z+a,\; g_{(a,\mu',\mu'')}(w)=\mu' w, \;g_{(a,\mu',\mu'')}
(\sqrt{-1})=\mu''\sqrt{-1}.$$ (In the case where $\sqrt{-1}$ lies in
${\mathbb F}_p$, we have ${\mathbb F}_p(z,w,\sqrt{-1})={\mathbb
F}_p(z,w)$, and we have an isomorphism
$${\mathbb F}_p \times \{\pm 1\}\stackrel \cong\to
H=\mathrm {Gal}({\mathbb F}_p(z,w)/{\mathbb F}_p(y))$$ which maps
$(a,\mu)\in {\mathbb F}_p\times\{\pm 1\}$ to the element
$g_{(a,\mu)}\in \mathrm {Gal}({\mathbb F}_p(z,w)/{\mathbb F}_p(y))$
defined by
$$g_{(a,\mu)}(z)=z+a,\; g_{(a,\mu)}(w)=\mu w.$$ The following argument works
for this case with slight modification. We leave to the reader to
treat this case.) Let $V$ be a one dimensional $\overline{\mathbb
Q}_\ell$-vector space with a basis $e_0$. Define an action of $H$ on
$V$ by
$$g_{(a,\mu',\mu'')}(e_0)=\psi(-2a){\chi}(\mu'^{-1})e_0.$$ Then $[2]_\ast
({\mathcal L}_{\psi}(2t)\otimes {\mathcal L}_{{\chi}})$ is just the
composition of $\mathrm {Ind}_H^G(V)$ with the canonical
homomorphism
$$\mathrm {Gal}(\overline {{\mathbb F}_p(t)}/{\mathbb F}_p(t))\to
\mathrm {Gal}({\mathbb F}_p(z,w,\sqrt{-1})/{\mathbb F}_p(t))=G.$$
Let $g$ be the element in $G=\mathrm {Gal}({\mathbb
F}_p(z,w,\sqrt{-1})/{\mathbb F}_p(t))$ defined by $$g(z)=-z,\;
g(w)=\sqrt{-1}w,\; g(\sqrt{-1})=\sqrt{-1}.$$ Then the image of $g$
in $G/H$ is a generator of the cyclic group $G/H$. So $G$ is
generated by $g_{(a,\mu',\mu'')}\in H$ $((a,\mu',\mu'')\in {\mathbb
F}_p\times \{\pm 1\}\times\{\pm 1\})$ and $g$. The space $\mathrm
{Ind}_H^G(V)$ has a basis $\{e_0, e_1\}$ with
\begin{eqnarray*}
g(e_0)&=&e_1, \\
g_{(a,\mu',\mu'')}(e_0)&=&
\psi(-2a){\chi}(\mu'^{-1})e_0,\\
g_{(a,\mu',\mu'')}(e_1)&=&\psi(2a){\chi}(\mu'^{-1}\mu''^{-1})e_1,\\
g(e_1)&=&g^2(e_0)=g_{(0,-1,1)}(e_0)=-e_0.
\end{eqnarray*}
Suppose $k=2r$ is even. $\mathrm {Sym}^k(\mathrm {Ind}_H^G(V))$ has
a basis
$$\{e_1^k, g(e_1^{k}), e_0e_1^{k-1}, g(e_0e_1^{k-1}),\ldots,
e_0^{r-1}e_1^{r+1}, g(e_0^{r-1}e_1^{r+1}), e_0^re_1^r\},$$ and for
each $i=0,1,\ldots, r$, we have
\begin{eqnarray*}
g_{(a,\mu',\mu'')}(e_0^ie_1^{k-i})&=&
\psi(-2ia){\chi}(\mu'^{-i})\psi(2(k-i)a){\chi}(\mu'^{-(k-i)}\mu''^{-(k-i)})e_0^ie_1^{k-i}\\
&=&
\psi(2(k-2i)a){\chi}(\mu'^{-k}){\chi}(\mu''^{-(k-i)})e_0^ie_1^{k-i}.
\end{eqnarray*}
Using the fact that $k$ is even and ${\chi}^2=1$, we get
$$g_{(a,\mu',\mu'')}(e_0^ie_1^{k-i})=\psi(2(k-2i)a){\chi}(\mu''^{i})e_0^ie_1^{k-i}.$$
In particular, we have
$$g_{(a,\mu',\mu'')}(e_0^re_1^r)={\chi}(\mu''^{r})e_0^re_1^r.$$
Moreover, we have
$$g(e_0^re_1^r)=e_1^r(g(e_1^r))=(-1)^re_0^re_1^r.$$
It follows that
$$\mathrm {Sym}^k([2]_\ast ({\mathcal L}_{\psi}(2t)\otimes
{\mathcal L}_{\chi}))\cong ({\mathcal L}_{{\chi}^r}\otimes {\mathcal
L}_{\theta_1^{r}})\oplus \left(\bigoplus_{i=0}^{r-1}
[2]_\ast({\mathcal L}_{\psi}(2(2i-k)t)\otimes {\mathcal
L}_{\theta_1^{i}})\right).
$$
We have
$$[2]_\ast({\mathcal
L}_{\psi}(2(2i-k)t)\otimes {\mathcal L}_{\theta_1^{i}})\cong
[2]_\ast({\mathcal L}_\psi(2(2i-k)t)\otimes [2]^\ast {\mathcal
L}_{\theta_1^{i}})\cong [2]_\ast{\mathcal L}_\psi(2(2i-k)t)\otimes
{\mathcal L}_{\theta_1^{i}}.$$ So we have
$$\mathrm {Sym}^k([2]_\ast ({\mathcal L}_{\psi}(2t)\otimes
{\mathcal L}_{\chi}))\cong ({\mathcal L}_{{\chi}^r}\otimes {\mathcal
L}_{\theta_1^{r}})\oplus \left(\bigoplus_{i=0}^{r-1}
[2]_\ast{\mathcal L}_\psi(2(2i-k)t)\otimes {\mathcal
L}_{\theta_1^{i}}\right).
$$
Suppose $n=2r+1$ is odd. $\mathrm {Sym}^k(\mathrm {Ind}_H^G(V))$ has
a basis
$$\{e_1^k, g(e_1^{k}), e_0e_1^{k-1}, g(e_0e_1^{k-1}),\ldots,
e_0^{r}e_1^{r+1}, g(e_0^{r}e_1^{r+1})\}.$$ Using the same
calculation as above, we get
$$\mathrm {Sym}^k([2]_\ast ({\mathcal L}_{\psi}(2t)\otimes
{\mathcal L}_{\chi}))\cong \bigoplus_{i=0}^{r} [2]_\ast({\mathcal
L}_{\psi}(2(2i-k)t)\otimes{\mathcal L}_{{\chi}})\otimes {\mathcal
L}_{\theta_1^{i+1}}.
$$
Lemma 2.4 follows by twisting the above expressions of $\mathrm
{Sym}^k([2]_\ast({\mathcal L}_{\psi}(2t)\otimes {\mathcal
L}_{{\chi}}))$ by ${\mathcal L}_{\theta_0^k}$.
\end{proof}

\begin{lemma} Assume $a\in {\mathbb F}_p$ is nonzero. We have the following identities.

(i) $\epsilon( {\mathbb P}^1_{(\infty)}, \overline{\mathbb Q}_\ell,
dt|_{{\mathbb P}^1_{(\infty)}})=\frac{1}{p^2}$.

(ii) $\epsilon( {\mathbb P}^1_{(\infty)}, \overline{\mathbb Q}_\ell,
dt^2|_{{\mathbb P}^1_{(\infty)}})=\frac{1}{p^3}$.

(iii) $\epsilon( {\mathbb P}^1_{(\infty)}, j_\ast{\mathcal
L}_\chi|_{{\mathbb P}^1_{(\infty)}}, dt|_{{\mathbb
P}^1_{(\infty)}})=-\frac{g(\chi,\psi)}{p^2}$.

(iv) $\epsilon( {\mathbb P}^1_{(\infty)}, j_\ast{\mathcal
L}_\chi|_{{\mathbb P}^1_{(\infty)}}, dt^2|_{{\mathbb
P}^1_{(\infty)}})=-\frac{g(\chi,\psi)}{p^3}\left(\frac{-2}{p}\right)$.

(v) $\epsilon( {\mathbb P}^1_{(\infty)}, j_\ast({\mathcal
L}_\psi(at)\otimes {\mathcal L}_\chi)|_{{\mathbb P}^1_{(\infty)}},
dt^2|_{{\mathbb P}^1_{(\infty)}})=\frac{1}{p^2}
\left(\frac{2a}{p}\right)$.

(vi) $\epsilon( {\mathbb P}^1_{(\infty)}, j_\ast{\mathcal
L}_\psi(at)|_{{\mathbb P}^1_{(\infty)}}, dt^2|_{{\mathbb
P}^1_{(\infty)}})=\frac{1}{p^2}$.

(vii) $\epsilon( {\mathbb P}^1_{(\infty)}, [2]_\ast\overline
{\mathbb Q}_\ell|_{{\mathbb P}^1_{(\infty)}}, dt|_{{\mathbb
P}^1_{(\infty)}})=-\frac{g(\chi,\psi)}{p^4}.$

(viii) $\epsilon( {\mathbb P}^1_{(\infty)}, j_\ast[2]_\ast({\mathcal
L}_\psi(at)\otimes {\mathcal L}_\chi)|_{{\mathbb P}^1_{(\infty)}},
dt|_{{\mathbb
P}^1_{(\infty)}})=-\frac{g(\chi,\psi)}{p^3}\left(\frac{2a}{p}\right)$.

(ix) $\epsilon( {\mathbb P}^1_{(\infty)}, j_\ast[2]_\ast{\mathcal
L}_\chi|_{{\mathbb P}^1_{(\infty)}}, dt|_{{\mathbb
P}^1_{(\infty)}})=\frac{g(\chi,\psi)^2}{p^4}\left(\frac{-2}{p}\right)$.

(x) $\epsilon( {\mathbb P}^1_{(\infty)}, j_\ast[2]_\ast {\mathcal
L}_\psi(at)|_{{\mathbb P}^1_{(\infty)}}, dt|_{{\mathbb
P}^1_{(\infty)}})=-\frac{g(\chi,\psi)}{p^3}$.
\end{lemma}

\begin{proof} Let $K_\infty$ be the completion of the field
$k(\eta_\infty)$, let ${\mathcal O}_\infty$ be the ring of integers
in $K_\infty$, and let $s=\frac{1}{t}$. Then $s$ is a uniformizer of
$K_\infty$. Denote the inclusion $\eta_\infty\to {\mathbb
P}^1_{(\infty)}$ also by $j$. Let $V$ be a $\overline{\mathbb
Q}_\ell$-sheaf of rank $1$ on $\eta_\infty$, and let
$\phi:K_\infty^\ast\to \overline{\mathbb Q}_\ell^\ast$ be the
character corresponding to $V$ via the reciprocity law. The Artin
conductor $a(\phi)$ of $\phi$ is defined to be the smallest integer
$m$ such that $\phi|_{1+s^m{\mathcal O}_\infty}=1$. For any nonzero
meromorphic differential $1$-form $\omega=fds$ on ${\mathbb
P}^1_{(\infty)}$, define the order $v_\infty(\omega)$ of $\omega$ to
be the valuation $v_\infty(f)$ of $f$. By \cite{L} 3.1.5.4 (v), we
have
\[\epsilon({\mathbb P}^1_{(\infty)}, j_\ast V, \omega)=
\left\{\begin {array} {ll}
\phi(s^{v_\infty(\omega)})p^{v_\infty(\omega)} & \hbox { if }
\phi|_{{\mathcal O}_\infty^\ast}=1, \\
\int_{s^{-(a(\phi)+v_\infty(\omega))}{\mathcal O}_\infty^\ast}
\phi^{-1}(z)\psi(\mathrm {Res}_\infty(z\omega))dz & \hbox{ if
}\phi|_{{\mathcal O}_\infty^\ast}\not=1,
\end{array}
\right.\] where $\mathrm {Res}_\infty$ denotes the residue of a
meromorphic $1$-form at $\infty$, and the integral is taken with
respect to the Haar measure $dz$ on $K_\infty$ normalized by
$\int_{{\mathcal O}_\infty}dz=1$.

Note that $dt=-\frac{ds}{s^2}$ has order $-2$ at $\infty$ and
$dt^2=-\frac{2ds}{s^3}$ has order $-3$. Applying the first case of
the above formula for the $\epsilon$-factor, we get (i) and (ii).

(iii) Taking $a=t=\frac{1}{s}$ and $b=z$ in the explicit reciprocity
law in \cite{S} XIV \S 3 Proposition 8, we see the character
$$\chi':K_\infty^\ast\to \overline {\mathbb Q}_l^\ast$$
corresponding to ${\mathcal L}_\chi$ is given by
$$\chi'(z)=\chi^{-1}\left(\left(\frac{\bar c}{p}\right)\right),$$
where $$c=(-1)^{-v_\infty(z)}\frac{z^{-1}}{s^{-v_\infty(z)}}$$ which
is a unit in ${\mathcal O}_\infty$,  $\bar c$ is the residue class
of $c$ in ${\mathcal O}_\infty/s{\mathcal O}_\infty\cong {\mathbb
F}_p$, and $\left(\frac{\bar c}{p}\right)$ is the Legendre symbol of
$\bar c$. Note that our formula for $c$ is the reciprocal of the
formula in \cite{S} because the reciprocity map in \cite{S} maps
uniformizers in $K$ to arithmetic Frobenius elements in $\mathrm
{Gal}(\overline K_\infty/K_\infty)^{\mathrm {ab}}$, whereas the
reciprocity map in \cite{L} maps uniformizers in $K$ to geometric
Frobenius elements. One can verify $a(\chi')=1$. For any $z\in
s{\mathcal O}_\infty^\ast$, write
$$z=s(r_0+r_1s+\cdots)$$ with $r_i\in {\mathbb F}_p$ and
$r_0\not=0$. We then have
\begin{eqnarray*}
\bar c&=&-r_0^{-1},\\
\mathrm{Res}_\infty(zdt)&=&-r_0.
\end{eqnarray*}
So we have
\begin{eqnarray*}
\epsilon( {\mathbb P}^1_{(\infty)}, j_\ast{\mathcal
L}_\chi|_{{\mathbb P}^1_{(\infty)}}, dt|_{{\mathbb
P}^1_{(\infty)}})&=&\int_{s{\mathcal O}_\infty^\ast}
\chi'^{-1}(z)\psi(\mathrm{Res}_\infty(zdt)) dz\\
&=& \int_{s{\mathcal O}_\infty^\ast}
\chi\left(\left(\frac{-r_0^{-1}}{p}\right)\right)\psi(-r_0)dz\\
&=&\int_{s{\mathcal O}_\infty^\ast}
\left(\frac{-r_0}{p}\right)\psi(-r_0)dz \\
&=&\sum_{r_0\in{\mathbb F}_p^\ast} \int_{r_0s(1+s{\mathcal
O}_\infty)}
\left(\frac{-r_0}{p}\right)\psi(-r_0)dz \\
&=& \sum_{r_0\in{\mathbb
F}_p^\ast}\left(\frac{-r_0}{p}\right)\psi(-r_0)\int_{r_0s(1+s{\mathcal
O}_\infty)}dz \\
&=& \frac{1}{p^2}\sum_{r_0\in{\mathbb
F}_p^\ast}\left(\frac{-r_0}{p}\right)\psi(-r_0)\\
&=&-\frac{g(\chi,\psi)}{p^2}.
\end{eqnarray*}

(iv) We can use the same method as in (iii), or use the formula
\cite{L} 3.1.5.5 to get
\begin{eqnarray*}
\epsilon( {\mathbb P}^1_{(\infty)}, j_\ast{\mathcal
L}_\chi|_{{\mathbb P}^1_{(\infty)}}, dt^2|_{{\mathbb
P}^1_{(\infty)}})&=& \epsilon( {\mathbb P}^1_{(\infty)},
j_\ast{\mathcal L}_\chi|_{{\mathbb P}^1_{(\infty)}}, 2tdt|_{{\mathbb
P}^1_{(\infty)}})\\
&=&\chi'\left(\frac{2}{s}\right)p^{v_\infty(\frac{2}{s})} \epsilon(
{\mathbb P}^1_{(\infty)}, j_\ast{\mathcal L}_\chi|_{{\mathbb
P}^1_{(\infty)}}, dt|_{{\mathbb P}^1_{(\infty)}})\\
&=&\left(\frac{-{2}}{p}\right)\cdot \frac{1}{p}\cdot \epsilon(
{\mathbb P}^1_{(\infty)}, j_\ast{\mathcal L}_\chi|_{{\mathbb
P}^1_{(\infty)}}, dt|_{{\mathbb P}^1_{(\infty)}})\\
&=&-\frac{g(\chi,\psi)}{p^3}\left(\frac{-2}{p}\right)
\end{eqnarray*}

(v) Taking $a$ to be $at=\frac{a}{s}$ and $b=z$ in the explicit
reciprocity law in \cite{S} XIV \S 5 Proposition 15, we see the
character
$$K_\infty^\ast\to \overline {\mathbb Q}_l^\ast$$
corresponding to ${\mathcal L}_\chi(at)$ is
$$z\mapsto \psi^{-1}\left(
-\mathrm
{Res}_\infty\left(\frac{a}{s}\cdot\frac{dz}{z}\right)\right).$$ (We
add the negative sign to the formula in \cite{S} since the
reciprocity map in \cite{S} is different from the one used in
\cite{L}.) So the character
$$\phi:K_\infty^\ast\to \overline {\mathbb Q}_l^\ast$$
corresponding to ${\mathcal L}_\chi(at)\otimes {\mathcal L}_\chi$ is
given by
$$\phi(z)= \psi^{-1}\left(
-\mathrm
{Res}_\infty\left(\frac{a}{s}\cdot\frac{dz}{z}\right)\right)
\chi^{-1}\left(\left(\frac{\bar c}{p}\right)\right),$$ where
$c=(-1)^{-v_\infty(z)}\frac{z^{-1}}{s^{-v_\infty(z)}}$. One can
verify $a(\phi)=2$. For any $z\in s{\mathcal O}_\infty^\ast$, write
$$z=s(r_0+r_1s+\cdots)$$ with $r_i\in {\mathbb F}_p$ and
$r_0\not=0$. We then have
\begin{eqnarray*}
\mathrm {Res}_\infty\left(\frac{a}{s}\cdot\frac{dz}{z}\right)&=&
\frac{ar_1}{r_0},\\
\bar c&=&-r_0^{-1},\\
\mathrm{Res}_\infty(zdt^2)&=&-2r_1.
\end{eqnarray*}
So we have
\begin{eqnarray*}
&&\epsilon( {\mathbb P}^1_{(\infty)}, j_\ast{\mathcal
L}_\chi|_{{\mathbb P}^1_{(\infty)}}, dt^2|_{{\mathbb
P}^1_{(\infty)}})\\&=&\int_{s{\mathcal O}_\infty^\ast}
\phi^{-1}(z)\psi(\mathrm{Res}_\infty(zdt^2)) dz\\
&=& \int_{s{\mathcal O}_\infty^\ast}
\psi\left(-\frac{ar_1}{r_0}\right)
\chi\left(\left(\frac{-r_0^{-1}}{p}\right)\right)\psi(-2r_1)dz\\
&=&\int_{s{\mathcal O}_\infty^\ast}
\left(\frac{-r_0}{p}\right)\psi\left(-r_1\left(\frac{a}{r_0}+2\right)\right)dz \\
&=&\sum_{r_0,r_1\in{\mathbb F}_p,r_0\not=0}
\int_{s(r_0+r_1s)(1+s^2{\mathcal O}_\infty)}
\left(\frac{-r_0}{p}\right)\psi\left(-r_1\left(\frac{a}{r_0}+2\right)\right)dz \\
&=& \sum_{r_0,r_1\in{\mathbb F}_p,r_0\not=0}
\left(\frac{-r_0}{p}\right)\psi\left(-r_1\left(\frac{a}{r_0}+2\right)\right)
\int_{s(r_0+r_1s)(1+s^2{\mathcal O}_\infty)}dz  \\
&=& \frac{1}{p^3}\sum_{r_0,r_1\in{\mathbb F}_p,r_0\not=0}
\left(\frac{-r_0}{p}\right)\psi\left(-r_1\left(\frac{a}{r_0}+2\right)\right)\\
&=&\frac{1}{p^3}\sum_{r_0\in{\mathbb F}_p^\ast}
\left(\frac{-r_0}{p}\right)\sum_{r_1\in{\mathbb
F}_p}\psi\left(-r_1\left(\frac{a}{r_0}+2\right)\right)\\
&=&\frac{1}{p^3}\cdot \left(\frac{-\frac{-a}{2}}{p}\right)\cdot p\\
&=& \frac{1}{p^2}\left(\frac{2a}{p}\right).
\end{eqnarray*}

We omit the proof of (vi), which is similar to the proof of (v).

(vii) We have $[2]_\ast\overline {\mathbb Q}_\ell\cong\overline
{\mathbb Q}_\ell\oplus j_\ast{\mathcal L}_\chi.$ So
$$\epsilon( {\mathbb P}^1_{(\infty)}, [2]_\ast\overline {\mathbb
Q}_\ell|_{{\mathbb P}^1_{(\infty)}}, dt|_{{\mathbb
P}^1_{(\infty)}})=\epsilon( {\mathbb P}^1_{(\infty)},\overline
{\mathbb Q}_\ell|_{{\mathbb P}^1_{(\infty)}}, dt|_{{\mathbb
P}^1_{(\infty)}})\epsilon( {\mathbb P}^1_{(\infty)}, j_\ast{\mathcal
L}_\chi|_{{\mathbb P}^1_{(\infty)}}, dt|_{{\mathbb
P}^1_{(\infty)}}).
$$ We then use (i) and (iii).

(viii) We can define $\epsilon$-factors for virtual sheaves on
${\mathbb P}^1_{(\infty)}$. By \cite{L} 3.1.5.4 (iv), we have
$$\epsilon( {\mathbb P}^1_{(\infty)}, [2]_\ast ([j_\ast ({\mathcal
L}_\psi(at)\otimes {\mathcal L}_\chi)]-[\overline {\mathbb Q}_\ell])
|_{{\mathbb P}^1_{(\infty)}}, dt|_{{\mathbb P}^1_{(\infty)}})=
\epsilon( {\mathbb P}^1_{(\infty)}, ([j_\ast ({\mathcal
L}_\psi(at)\otimes {\mathcal L}_\chi)]-[\overline {\mathbb Q}_\ell])
|_{{\mathbb P}^1_{(\infty)}}, dt^2|_{{\mathbb P}^1_{(\infty)}}).$$
Hence
\begin{eqnarray*}
&&\epsilon( {\mathbb P}^1_{(\infty)}, j_\ast[2]_\ast ({\mathcal
L}_\psi(at)\otimes {\mathcal L}_\chi)
|_{{\mathbb P}^1_{(\infty)}}, dt|_{{\mathbb P}^1_{(\infty)}})\\
&=&\frac{\epsilon( {\mathbb P}^1_{(\infty)}, j_\ast ({\mathcal
L}_\psi(at)\otimes {\mathcal L}_\chi)|_{{\mathbb P}^1_{(\infty)}},
dt^2|_{{\mathbb P}^1_{(\infty)}})}{\epsilon( {\mathbb
P}^1_{(\infty)}, \overline {\mathbb Q}_\ell , dt^2|_{{\mathbb
P}^1_{(\infty)}})} \epsilon( {\mathbb P}^1_{(\infty)},
[2]_\ast\overline{\mathbb Q}_\ell, dt|_{{\mathbb P}^1_{(\infty)}}).
\end{eqnarray*}
We then apply the formulas (ii), (v), and (vii).

We omit the proof of (ix) and (x), which is similar to the
proof of (viii).
\end{proof}

\begin{lemma} We have
\begin{eqnarray*} &&\epsilon( {\mathbb P}^1_{(\infty)},
j_\ast({\mathcal L}_{\chi^r}\otimes {\mathcal
L}_{\theta_0^{2r}\theta_1^r})
|_{{\mathbb P}^1_{(\infty)}}, dt|_{{\mathbb P}^1_{(\infty)}})\\
&=&\left\{
\begin {array}{ll}
\frac{g(\chi,\psi)^{-4r}}{p^2}& \hbox { if } r \hbox
{ is even},\\
-\frac{g(\chi,\psi)^{-2r+1}}{p^2}\left(\frac{-1}{p}\right)& \hbox {
if } r \hbox { is odd},
\end{array}
\right.\\ &&\epsilon( {\mathbb P}^1_{(\infty)}, j_\ast
([2]_\ast{\mathcal L}_{\psi}((4i-4r)t)\otimes {\mathcal
L}_{\theta_0^{2r}\theta_1^i})
|_{{\mathbb P}^1_{(\infty)}}, dt|_{{\mathbb P}^1_{(\infty)}})\\
&=&\left\{
\begin {array}{ll}
-\frac{g(\chi,\psi)^{-6r+1}}{p^4}\left(\frac{-1}{p}\right)^i& \hbox { if } p|i-r,\\
-\frac{g(\chi,\psi)^{-2r+1}}{p^3}\left(\frac{-1}{p}\right)^i & \hbox
{ if } p\not | i-r,
\end{array}
\right.\\
&&\epsilon( {\mathbb P}^1_{(\infty)}, j_\ast([2]_\ast{\mathcal
L}_\psi((4i-4r-2)t)\otimes {\mathcal L}_{\chi})\otimes {\mathcal
L}_{\theta_0^{2r+1}\theta_1^{i+1}}) |_{{\mathbb P}^1_{(\infty)}},
dt|_{{\mathbb P}^1_{(\infty)}}) \\
&=&\left\{
\begin {array}{ll}
\frac{g(\chi,\psi)^{-4r}}{p^4}\left(\frac{-2}{p}\right)& \hbox { if } p|2i-2r-1,\\
-\frac{g(\chi,\psi)^{-2r}}{p^3}\left(\frac{(-1)^{i+1}(2i-2r-1)}{p}\right)&
\hbox { if } p\not |2i-2r-1.
\end{array}
\right.
\end{eqnarray*}
\end{lemma}

\begin{proof} Let $F_\infty$ be the geometric Frobenius element at $\infty$. We have
$$\theta_0(F_\infty)=g(\chi,\psi), \;
\theta_1(F_\infty)=\left(\frac{-1}{p}\right).$$ Using the notation
\cite{L} 3.1.5.1, we have
\begin{eqnarray*}
a({\mathbb P}^1_{(\infty)},\overline {\mathbb Q}_\ell,dt|_{{\mathbb
P}^1_{(\infty)}})&=& -2,\\
a({\mathbb P}^1_{(\infty)},j_\ast {\mathcal L}_\chi|_{{\mathbb
P}^1_{(\infty)}},dt|_{{\mathbb
P}^1_{(\infty)}})&=& -1,\\
a({\mathbb P}^1_{(\infty)},[2]_\ast \overline {\mathbb
Q}_\ell,dt|_{{\mathbb
P}^1_{(\infty)}})&=& -3,\\
a({\mathbb P}^1_{(\infty)},j_\ast[2]_\ast {\mathcal
L}_\psi(at)|_{{\mathbb P}^1_{(\infty)}},dt|_{{\mathbb
P}^1_{(\infty)}})&=& -1,\; (a\in {\mathbb F}_p^\ast)\\
a({\mathbb P}^1_{(\infty)},j_\ast[2]_\ast {\mathcal
L}_\chi|_{{\mathbb P}^1_{(\infty)}},dt|_{{\mathbb
P}^1_{(\infty)}})&=& -2,\\
a({\mathbb P}^1_{(\infty)},j_\ast[2]_\ast ({\mathcal
L}_\psi(at)\otimes{\mathcal L}_\chi)|_{{\mathbb
P}^1_{(\infty)}},dt|_{{\mathbb P}^1_{(\infty)}})&=& -1,\; (a\in
{\mathbb F}_p^\ast).
\end{eqnarray*}
So by \cite{L} 3.1.5.6, we have
\begin{eqnarray*} &&\epsilon( {\mathbb P}^1_{(\infty)},
j_\ast({\mathcal L}_{\chi^r}\otimes {\mathcal
L}_{\theta_0^{2r}\theta_1^r})
|_{{\mathbb P}^1_{(\infty)}}, dt|_{{\mathbb P}^1_{(\infty)}})\\
&=&\left\{
\begin {array}{ll}
((\theta_0^{2r}\theta_1^r)(F_\infty))^{-2}\epsilon( {\mathbb
P}^1_{(\infty)}, \overline{\mathbb Q}_\ell, dt|_{{\mathbb
P}^1_{(\infty)}})& \hbox { if } r \hbox
{ is even},\\
((\theta_0^{2r}\theta_1^r)(F_\infty))^{-1}\epsilon({\mathbb
P}^1_{(\infty)}, j_\ast{\mathcal L}_{\chi}|_{{\mathbb
P}^1_{(\infty)}}, dt|_{{\mathbb P}^1_{(\infty)}})& \hbox { if } r
\hbox { is odd},
\end{array}
\right.\\
&&\epsilon( {\mathbb P}^1_{(\infty)}, j_\ast([2]_\ast{\mathcal
L}_{\psi}((4i-4r)t)\otimes {\mathcal L}_{\theta_0^{2r}\theta_1^i})
|_{{\mathbb P}^1_{(\infty)}}, dt|_{{\mathbb P}^1_{(\infty)}})\\
&=&\left\{
\begin {array}{ll}
((\theta_0^{2r}\theta_1^i)(F_\infty))^{-3}\epsilon( {\mathbb
P}^1_{(\infty)}, [2]_\ast\overline{\mathbb Q}_\ell,
dt|_{{\mathbb P}^1_{(\infty)}})& \hbox { if } p|i-r,\\
((\theta_0^{2r}\theta_1^i)(F_\infty))^{-1}\epsilon( {\mathbb
P}^1_{(\infty)}, j_\ast [2]_\ast{\mathcal
L}_{\psi}((4i-4r)t)|_{{\mathbb P}^1_{(\infty)}}, dt|_{{\mathbb
P}^1_{(\infty)}})& \hbox { if } p\not | i-r,
\end{array}
\right.\\
&&\epsilon( {\mathbb P}^1_{(\infty)}, j_\ast ([2]_\ast({\mathcal
L}_\psi((4i-4r-2)t)\otimes {\mathcal L}_{\chi})\otimes {\mathcal
L}_{\theta_0^{2r+1}\theta_1^{i+1}}) |_{{\mathbb P}^1_{(\infty)}},
dt|_{{\mathbb P}^1_{(\infty)}}) \\
&=&\left\{
\begin {array}{ll}
((\theta_0^{2r+1}\theta_1^{i+1})(F_\infty))^{-2}\epsilon( {\mathbb
P}^1_{(\infty)}, j_\ast [2]_\ast{\mathcal L}_\chi|_{{\mathbb
P}^1_{(\infty)}},
dt|_{{\mathbb P}^1_{(\infty)}})& \hbox { if } p|2i-2r-1,\\
((\theta_0^{2r+1}\theta_1^{i+1})(F_\infty))^{-1}\epsilon( {\mathbb
P}^1_{(\infty)}, j_\ast [2]_\ast({\mathcal
L}_\psi((4i-4r-2)t)\otimes {\mathcal L}_{\chi})|_{{\mathbb
P}^1_{(\infty)}}, dt|_{{\mathbb P}^1_{(\infty)}}) & \hbox { if }
p\not |2i-2r-1.
\end{array}
\right.
\end{eqnarray*}
We then apply the formulas in Lemma 2.5.
\end{proof}

The following is Proposition 0.4 in the introduction.

\begin{proposition}
$\epsilon( {\mathbb P}^1_{(\infty)}, j_\ast(\mathrm {Sym}^k(\mathrm
{Kl}_2))|_{{\mathbb P}^1_{(\infty)}}, dt|_{{\mathbb
P}^1_{(\infty)}})$ equals
$$p^{-(k+1)(\frac{k+8}{4}+[\frac{k}{2p}])}$$
if $k=2r$ for an even $r$,
$$p^{-(k+1)(\frac{k+6}{4}+[\frac{k}{2p}])}$$
if $k=2r$ for an odd $r$, and
$$(-1)^{\frac{k+1}{2}+[\frac{k}{p}]-[\frac{k}{2p}]}p^{-\frac{k+1}{2}(\frac{k+5}{2}+[\frac{k}{p}]-[\frac{k}{2p}])}
\left( \frac{-2}{p}\right)^{[\frac{k}{p}]-[\frac{k}{2p}]}
\prod_{j\in\{0,1,\ldots,[\frac{k}{2}]\},\; p\not\; | 2j+1}\left(
\frac{(-1)^{j}(2j+1)}{p}\right)
$$
if $k=2r+1$.
\end{proposition}

\begin{proof} By Lemmas 2.4 and 2.6,
$\epsilon( {\mathbb P}^1_{(\infty)}, j_\ast(\mathrm {Sym}^k(\mathrm
{Kl}_2))|_{{\mathbb P}^1_{(\infty)}}, dt|_{{\mathbb
P}^1_{(\infty)}})$ equals
$$\frac{g(\chi,\psi)^{-4r}}{p^2}\prod_{i\in\{0,\ldots,r-1\},\;p|i-r}
\left(-\frac{g(\chi,\psi)^{-6r+1}}{p^4}\left(\frac{-1}{p}\right)^i\right)
\prod_{i\in\{0,\ldots,r-1\},\;p\not\;|i-r}\left(-\frac{g(\chi,\psi)^{-2r+1}}{p^3}\left(\frac{-1}{p}\right)^i
\right)$$ if $k=2r$ for an even $r$,
$$-\frac{g(\chi,\psi)^{-2r+1}}{p^2}\left(\frac{-1}{p}\right)
\prod_{i\in\{0,\ldots,r-1\},\;p|i-r}
\left(-\frac{g(\chi,\psi)^{-6r+1}}{p^4}\left(\frac{-1}{p}\right)^i\right)
\prod_{i\in\{0,\ldots,r-1\},\;p\not\;|i-r}\left(
-\frac{g(\chi,\psi)^{-2r+1}}{p^3}\left(\frac{-1}{p}\right)^i\right)
$$ if $k=2r$ for an odd $r$, and
$$\prod_{i\in\{0,\ldots,r\},\;p|2i-2r-1}
\left(\frac{g(\chi,\psi)^{-4r}}{p^4}\left(\frac{-2}{p}\right)\right)
\prod_{i\in\{0,\ldots,r\},\;p\not\;|2i-2r-1}
\left(-\frac{g(\chi,\psi)^{-2r}}{p^3}\left(\frac{(-1)^{i+1}(2i-2r-1)}{p}\right)\right)$$
if $k=2r+1$. Let's simplify the above expressions. Recall that
$g(\chi,\psi)^2=p\left(\frac{-1}{p}\right).$ If $k=2r$ with $r$
even, we have
\begin{eqnarray*}
&&\epsilon( {\mathbb P}^1_{(\infty)}, j_\ast(\mathrm {Sym}^k(\mathrm
{Kl}_2))|_{{\mathbb P}^1_{(\infty)}}, dt|_{{\mathbb
P}^1_{(\infty)}})\\
&=&\frac{g(\chi,\psi)^{-4r}}{p^2}\prod_{i\in\{0,\ldots,r-1\},\;p|i-r}
\left(-\frac{g(\chi,\psi)^{-6r+1}}{p^4}\left(\frac{-1}{p}\right)^i\right)
\prod_{i\in\{0,\ldots,r-1\},\;p\not\;|i-r}\left(-\frac{g(\chi,\psi)^{-2r+1}}{p^3}\left(\frac{-1}{p}\right)^i
\right)\\
&=&\frac{g(\chi,\psi)^{-4r}}{p^2}
\prod_{i\in\{0,\ldots,r-1\},\;p|i-r} \frac{g(\chi,\psi)^{-4r}}{p}
\prod_{i\in\{0,\ldots,r-1\}}\left(-\frac{g(\chi,\psi)^{-2r+1}}{p^3}\left(\frac{-1}{p}\right)^i
\right)\\
&=& \frac{g(\chi,\psi)^{-4r}}{p^2} \left(
\frac{g(\chi,\psi)^{-4r}}{p}\right)^{[\frac{r}{p}]}
\left(-\frac{g(\chi,\psi)^{-2r+1}}{p^3}\right)^{r}\left( \frac{-1}{p}\right)^{\frac{r(r-1)}{2}}\\
&=&\frac{p^{-2r}}{p^2}\left(\frac{p^{-2r}}{p}\right)^{[\frac{r}{p}]}
\frac{\left(p\left(\frac{-1}{p}\right)\right)^{\frac{r(-2r+1)}{2}}}{p^{3r}}\left(\frac{-1}{p}\right)
^{\frac{r(r-1)}{2}}\\
&=&
p^{-r^2-\frac{9}{2}r-2-(2r+1)[\frac{r}{p}]}\left(\frac{-1}{p}\right)^{-\frac{r^2}{2}}
\\
&=&p^{-(k+1)(\frac{k+8}{4}+[\frac{k}{2p}])}.
\end{eqnarray*}
If $k=2r$ with $r$ odd, we have
\begin{eqnarray*}
&&\epsilon( {\mathbb P}^1_{(\infty)}, j_\ast(\mathrm {Sym}^k(\mathrm
{Kl}_2))|_{{\mathbb P}^1_{(\infty)}}, dt|_{{\mathbb
P}^1_{(\infty)}})\\
&=& -\frac{g(\chi,\psi)^{-2r+1}}{p^2}\left(\frac{-1}{p}\right)
\prod_{i\in\{0,\ldots,r-1\},\;p|i-r}
\left(-\frac{g(\chi,\psi)^{-6r+1}}{p^4}\left(\frac{-1}{p}\right)^i\right)
\\&&\qquad\times \prod_{i\in\{0,\ldots,r-1\},\;p\not\;|i-r}\left(
-\frac{g(\chi,\psi)^{-2r+1}}{p^3}\left(\frac{-1}{p}\right)^i\right)\\
&=& -\frac{g(\chi,\psi)^{-2r+1}}{p^2}\left(\frac{-1}{p}\right)
\prod_{i\in\{0,\ldots,r-1\},\;p|i-r}
\left(\frac{g(\chi,\psi)^{-4r}}{p}\right)
\prod_{i\in\{0,\ldots,r-1\}}\left(
-\frac{g(\chi,\psi)^{-2r+1}}{p^3}\left(\frac{-1}{p}\right)^i\right)\\
&=& -\frac{g(\chi,\psi)^{-2r+1}}{p^2}\left(\frac{-1}{p}\right)
\left(\frac{g(\chi,\psi)^{-4r}}{p}\right)^{[\frac{r}{p}]} \left(
-\frac{g(\chi,\psi)^{-2r+1}}{p^3}\right)^r \left(\frac{-1}{p}\right)^{\frac{r(r-1)}{2}}\\
&=&
\frac{g(\chi,\psi)^{(-2r+1)(r+1)}}{p^{3r+2}}\left(\frac{-1}{p}\right)^{1+\frac{r(r-1)}{2}}
\left(\frac{g(\chi,\psi)^{-4r}}{p}\right)^{[\frac{r}{p}]} \\
&=&
\frac{\left(p\left(\frac{-1}{p}\right)\right)^{\frac{(-2r+1)(r+1)}{2}}}{p^{3r+2}}
\left(\frac{-1}{p}\right)^{1+\frac{r(r-1)}{2}}
\left(\frac{p^{-2r}}{p}\right)^{[\frac{r}{p}]} \\
&=&p^{-r^2-\frac{7}{2}r-\frac{3}{2}-(2r+1)[\frac{r}{p}]}\left(\frac{-1}{p}\right)^{-\frac{(r-1)(r+3)}{2}}\\
&=&p^{-(k+1)(\frac{k+6}{4}+[\frac{k}{2p}])}.
\end{eqnarray*}
If $k=2r+1$ is odd, we have
\begin{eqnarray*}
&&\epsilon( {\mathbb P}^1_{(\infty)}, j_\ast(\mathrm {Sym}^k(\mathrm
{Kl}_2))|_{{\mathbb P}^1_{(\infty)}}, dt|_{{\mathbb
P}^1_{(\infty)}})\\
&=&\prod_{i\in\{0,\ldots,r\},\;p|2i-2r-1}
\left(\frac{g(\chi,\psi)^{-4r}}{p^4}\left(\frac{-2}{p}\right)\right)
\prod_{i\in\{0,\ldots,r\},\;p\not\;|2i-2r-1}
\left(-\frac{g(\chi,\psi)^{-2r}}{p^3}\left(\frac{(-1)^{i+1}(2i-2r-1)}{p}\right)\right)\\
&=&\left(-\frac{g(\chi,\psi)^{-2r}}{p^3}\right)^{r+1}
\left(-\frac{g(\chi,\psi)^{-2r}}{p}\left(\frac{-2}{p}\right)\right)
^{[\frac{k}{p}]-[\frac{k}{2p}]}
\prod_{i\in\{0,\ldots,r\},\;p\not\;|2i-2r-1}
\left(\frac{(-1)^{i+1}(2i-2r-1)}{p}\right)\\
&=& \left( -p^{-3-r}\left( \frac{-1}{p}\right)^{-r}\right)^{r+1}
\left( -p^{-1-r}\left( \frac{-1}{p}\right)^{-r} \left(
\frac{-2}{p}\right)\right)^{[\frac{k}{p}]-[\frac{k}{2p}]}\prod_{j\in\{0,\ldots,r\},\;
p\not\; | 2j+1}\left( \frac{(-1)^{r-j}(2j+1)}{p}\right)\\
&=&(-1)^{r+1+[\frac{k}{p}]-[\frac{k}{2p}]}
p^{-(r+1)(r+3+[\frac{k}{p}]-[\frac{k}{2p}])}\left(
\frac{-2}{p}\right)^{[\frac{k}{p}]-[\frac{k}{2p}]}\left(
\frac{-1}{p}\right)^{-r(r+1)-r([\frac{k}{p}]-[\frac{k}{2p}])}\\
&&\qquad\times \prod_{j\in\{0,\ldots,r\}, p\not\; | 2j+1}\left(
\frac{(-1)^{r-j}(2j+1)}{p}\right)\\
&=&
(-1)^{\frac{k+1}{2}+[\frac{k}{p}]-[\frac{k}{2p}]}p^{-\frac{k+1}{2}(\frac{k+5}{2}+[\frac{k}{p}]-[\frac{k}{2p}])}
\left( \frac{-2}{p}\right)^{[\frac{k}{p}]-[\frac{k}{2p}]}
\prod_{j\in\{0,\ldots,[\frac{k}{2}]\},\; p\not\; | 2j+1}\left(
\frac{(-1)^{j}(2j+1)}{p}\right).
\end{eqnarray*}
\end{proof}

\end{document}